\documentclass[a4,amstex,11pt]{article}
\usepackage{amsmath,wrapfig}
\usepackage{amssymb,bm,epic}
\usepackage[mathscr]{eucal}
\usepackage{amsthm, amsfonts, latexsym,enumerate, algorithm, algpseudocode}
\usepackage{graphicx}
\usepackage{setspace}
\makeatletter
\def\fixedfigure{\def\@captype{figure}}
\def\fixedtable{\def\@captype{table}}
\numberwithin{equation}{section}
\theoremstyle{definition}

\newtheorem{definition}{Definition}[section]
\newtheorem{theorem}[definition]{Theorem}
\newtheorem{proposition}[definition]{Proposition}
\newtheorem{lemma}[definition]{Lemma}
\newtheorem{corollary}[definition]{Corollary}
\newtheorem{remark}[definition]{Remark}

\algnewcommand\algorithmicinput{{\bfseries Input:}}
\algnewcommand\AlgInput{\item[\algorithmicinput]}
\algnewcommand\algorithmicoutput{{\bfseries Output:}}
\algnewcommand\AlgOutput{\item[\algorithmicoutput]}
\title{Algorithm for the CSR expansion of max-plus matrices using the characteristic polynomial 
\thanks{This work is supported by JSPS KAKENHI No.22K13964.}
}

\author{
Yuki Nishida\footnote{Department of Science, Technology and Informatics, Kyoto Prefectural University, Kyoto 606-8522, Japan.
Email: y-nishida@kpu.ac.jp
}
{}\footnotemark[0]
}

\date{}
\begin{document}
\maketitle

% =====================================

%\renewcommand{\abstractname}{\large \textbf{Abstract}}
\begin{abstract}
Max-plus algebra is a semiring with addition $a\oplus b = \max(a,b)$ and multiplication $a\otimes b = a+b$. It is applied in cases, such as combinatorial optimization and discrete event systems. We consider the power of max-plus square matrices, which is equivalent to obtaining the all-pair maximum weight paths with a fixed length in the corresponding weighted digraph. Each $n$-by-$n$ matrix admits the CSR expansion that decomposes the matrix into a sum of at most $n$ periodic terms after $O(n^{2})$ times of powers. In this study, we propose an $O(n(m+n \log n))$ time algorithm for the CSR expansion, where $m$ is the number of nonzero entries in the matrix, which improves the $O(n^{4} \log n)$ algorithm known for this problem. Our algorithm is based on finding the roots of the characteristic polynomial of the max-plus matrix. These roots play a similar role to the eigenvalues of the matrix, and become the growth rates of the terms in the CSR expansion.
\end{abstract}

% =====================================

{\it Keywords}: max-plus algebra, tropical semiring, CSR expansion, Nachtigall expansion, algebraic eigenvalue, shortest path problem

{\it 2020MSC}: 15A80, 93C65

\section{Introduction}

Max-plus algebra $\mathbb{R}_{\max} = \mathbb{R}\cup \{\varepsilon \}$, where $\varepsilon := -\infty$, is a semiring with addition $\oplus$ and multiplication $\otimes$ defined by   
\begin{align*}
	a \oplus b = \max (a,b), \quad a \otimes b = a+b 
\end{align*}
for $a,b \in \mathbb{R}_{\max}$. Matrix arithmetic over max-plus algebra is also defined as in conventional linear algebra, for example, 
\begin{align*}
	[A \oplus B]_{ij} &= [A]_{ij} \oplus [B]_{ij}, \\
    [A \otimes B]_{ij} &= \bigoplus_{k} [A]_{ik} \otimes [B]_{kj}, \\
    [c \otimes A]_{ij} &= c \otimes [A]_{ij}
\end{align*}
for compatible size of matrices $A$ and $B$, and scalar $c$, where $[A]_{ij}$ denotes the $(i,j)$ entry of $A$.
Max-plus algebra has been applied to many problems, such as combinatorial optimization~\cite{Butkovic2003} and discrete event systems~\cite{Baccelli1992,Schutter2020,Komenda2018}. It appears in many real-world problems, for example, steelworks~\cite{Green1960}, train timetable~\cite{Vries1998, Farhi2017, Goverde1998}, operation in emergency call center~\cite{Allamigeon2015,Allamigeon2021}, multiprocessor interactive system~\cite{Butkovic2009}, etc. Furthermore, max-plus algebra is used as a language for algebraic geometry over valuated fields, called tropical geometry~\cite{Joswig2021, Maclagan2015}. One of the advantages of considering max-plus algebra is that there is a one-to-one correspondence between max-plus square matrices and weighted digraphs. Therefore, some algebraic concepts are related to network optimization problems. For example, the maximum eigenvalue of a matrix is identical to the maximum mean weight of circuits in the corresponding graph~\cite{Green1962}. 

Consider a discrete event system defined by
\begin{align*}
    \bm{x}(t) = A \otimes \bm{x}(t-1), \quad t=0,1,\dots,
\end{align*}
where $A \in \mathbb{R}_{\max}^{n \times n}$ and $\bm{x}(0) \in \mathbb{R}_{\max}^{n}$. The state vector $\bm{x}(t)$ is represented by $\bm{x}(t) = A^{\otimes t} \otimes \bm{x}(0)$, where $A^{\otimes t}$ is the $t$th power of $A$ with respect to max-plus matrix product. Thus, computation of the matrix power $A^{\otimes t}$ is important. As a fundamental result, it was shown in~\cite{Cohen1985} that any irreducible matrix is ultimately periodic, that is, there exist positive integers $T(A)$ and $\sigma$, and a real number $\lambda$, such that
\begin{align*}
    A^{\otimes (t+\sigma)} = \lambda^{\otimes \sigma} \otimes A^{\otimes t}
\end{align*}
for all $t \geq T(A)$. The minimum integer $T(A)$ is called the transient of $A$. Although it completely determines the asymptotic behavior of the system, the transient $T(A)$ depends on the entries of $A$, not only on the dimension $n$. Upper bounds for $T(A)$ were proposed in~\cite{CB2013,CB2012,Hartmann1999}. In addition, it was proved that the critical rows or columns of the matrix becomes periodic after $2n^{2}$th powers~\cite{Sergeev2009-2}. To obtain a transient bound that is a polynomial of $n$, Nachtigall~\cite{Nachtigall1997} introduced an expansion of $A^{\otimes t}$ of the form 
\begin{align*}
    A^{\otimes t} = A_{1}^{\otimes t} \oplus A_{2}^{\otimes t} \oplus \cdots \oplus A_{r}^{\otimes t},
\end{align*}
where every $A_{k}$ is periodic if $t \geq 3n^{2}$, that is, there exist $\sigma_{k}$ and $\lambda_{k}$ such that $A_{k}^{\otimes (t+\sigma_{k})} = \lambda_{k}^{\otimes \sigma_{k}} \otimes A^{\otimes t}$ for $t \geq 3n^{2}$. The number $r$ of the terms in the expansion is at most $n$. The computational complexity of this expansion is $O(n^{6})$ in Nachtigall~\cite{Nachtigall1997}, and it was improved to $O(n^{5})$ by Moln\'{a}rov\'{a}~\cite{Molnarova2003}. In Sergeev and Schneider~\cite{Sergeev2012}, the Nachtigall expansion was reformulated as the CSR expansion:
\begin{align*}
    A^{\otimes t} = \bigoplus_{k} \lambda_{k}^{\otimes t} \otimes C_{k} \otimes S_{k}^{\otimes t} \otimes R_{k}, \quad t \geq 3n^{2}.
\end{align*}
Here, $C_{k}$ and $R_{k}$ are extracted from the Kleene star of a certain matrix, and $S_{k}$ is a permutation matrix (in max-plus algebraic sense). An $O(n^{4}\log n)$ algorithm for the CSR expansion is derived in~\cite{Sergeev2012}. In terms of the CSR expansion, several types of bounds for the transients of max-plus matrices have been proposed in~\cite{KCP2021,Merlet2014,Merlet2014-2}. An application of the CSR expansion to a tropical public key exchange protocol was developed in~\cite{Muanalifah2022}. Recently, the CSR expansion was extended to a sequence of inhomogeneous matrices~\cite{KCP2022}. 

In this study, we propose an $O(n(m+n \log n))$ algorithm for the CSR expansion of $A \in \mathbb{R}_{\max}^{n \times n}$, where $m$ is the number of finite (non-$\varepsilon$) entries in $A$. As a new scheme for the CSR expansion, we exploit the characteristic polynomial $\chi_{A}(t)$ of $A$~\cite{Green1983}. We show that each $\lambda_{k}$ in the CSR expansion is a root of $\chi_{A}(t)$. Note that the roots of $\chi_{A}(t)$ are called algebraic eigenvalues in~\cite{Akian2005}, and have properties similar to the usual eigenvalues. For example, every algebraic eigenvalue has a corresponding vector, which is called an algebraic eigenvector~\cite{Nishida2020}. These roots are also considered in supertropical algebra~\cite{Izhakian2011,Niv2017}, which is an extension of max-plus algebra. Gassner and Klinz~\cite{Gassner2009} proposed an $O(n(m+n \log n))$ algorithm to obtain all roots of $\chi_{A}(t)$ by reducing the problem to the parametric shortest path tree problem~\cite{Young1991}. This contributes in reducing the computational complexity of the CSR expansion because obtaining all $\lambda_{k}$ in the CSR expansion takes $O(mn^{2})$ time if we repeatedly apply the Karp Algorithm for the maximum circuit mean problem~\cite{Karp1978}. We focus on the maximum weight path problem appearing in computing the CSR expansion, which is equivalent to the shortest path problem by reversing the signs of the arc weight. The problem can be solved efficiently by the Dijkstra Algorithm~\cite{Dijkstra1959} if the weights of all arcs are non-positive. The visualization introduced in~\cite{Sergeev2009} is a transformation of a max-plus matrix into a non-positive one equivalent to it. We propose an efficient method to visualize all concerning matrices by using the visualization of smaller ones. In particular, we apply the Dijkstra Algorithm successively, making the dimension of matrix to be large. Furthermore, we consider the maximum weight path problem modulo an integer $\ell$, and then solve this problem by extending the graph to $\ell$ copies of it. Subsequently, we show that the total computational complexity is bounded by $O(n(m+n \log n))$ from above.

The rest of this paper is organized as follows. In Section 2, we present the notation and preliminary results on max-plus linear algebra, especially the max-plus eigenvalue problem in terms of weighted digraph. In Section 3, we describe the CSR expansion introduced in~\cite{Sergeev2012}. In Section 4, we develop an $O(n(m+n \log n))$ algorithm for the CSR expansion. Our method is summarized into three points, and Sections 4.1, 4.2 and 4.3 are devoted to each of them. An example demonstrating the algorithm is presented in Section 4.4. In Section 5, we give concluding remarks.

\section{Preliminaries on max-plus algebra}

Max-plus algebra is the set $\mathbb{R}_{\max} = \mathbb{R}\cup \{\varepsilon \}$, where $\varepsilon := -\infty$, with two operations $\oplus$ and $\otimes$, defined by   
\begin{align*}
	a \oplus b = \max (a,b), \quad a \otimes b = a+b
\end{align*}
for $a,b \in \mathbb{R}_{\max}$. Considering $\oplus$ and $\otimes$ as addition and multiplication, respectively, max-plus algebra is a semiring. Here, $\varepsilon$ is the identity element for addition, and $e := 0$ is the identity element for multiplication. 

Let $\mathbb{R}_{\max}^n$ and $\mathbb{R}_{\max}^{m \times n }$ be the set of $n$-dimensional max-plus column vectors and the set of $m \times n$ max-plus matrices, respectively. The operations $\oplus$ and $\otimes$ are extended to max-plus vectors and matrices as in conventional linear algebra. For $A, B \in \mathbb{R}_{\max}^{m \times n}$, the matrix sum $A \oplus B \in \mathbb{R}_{\max}^{m \times n}$ is defined by
\begin{align*}
    [A \oplus B]_{ij} = [A]_{ij} \oplus [B]_{ij},
\end{align*}
where $[A]_{ij}$ indicates the $(i,j)$ entry of $A$. For $A \in \mathbb{R}_{\max}^{\ell \times m}$ and $B \in \mathbb{R}_{\max}^{m \times n}$, the matrix product $A \otimes B \in \mathbb{R}_{\max}^{\ell \times n}$ is defined by
\begin{align*}
    [A \otimes B]_{ij} = \bigoplus_{k=1}^{m} [A]_{ik} \otimes [B]_{kj}.
\end{align*}
For $A \in \mathbb{R}_{\max}^{m \times n}$ and $c \in \mathbb{R}$, the scalar multiplication of $A$ by $c$ is defined by
\begin{align*}
    [c \otimes A]_{ij} = c \otimes [A]_{ij}.
\end{align*}
The max-plus zero vector or the zero matrix is denoted by $\mathcal{E}$, and the max-plus unit matrix of order $n$ is denoted by $I_{n}$. A max-plus diagonal matrix $\mathrm{diag}(\bm{d})$ defined by a vector $\bm{d} = (d_{1},d_{2},\dots,d_{n})^{\top} \in \mathbb{R}^{n}$ is the matrix whose diagonal entries are $d_{1},d_{2},\dots, d_{n}$ and off-diagonal entries are $\varepsilon$. Its inverse is expressed by $\mathrm{diag}(-\bm{d})$. For a matrix $A = (a_{ij}) \in\mathbb{R}_{\max}^{n \times n}$, we define the determinant of $A$ by 
\begin{align*}
	\det A = \bigoplus_{\pi \in S_n} \bigotimes_{i=1}^n a_{i \pi(i)},
\end{align*}
where $S_{n}$ is the symmetric group of order $n$. Computation of the max-plus determinant is equivalent to solving the maximum weight assignment problem in a bipartite graph.
\subsection{Max-plus matrices and graphs}
For a matrix $A = (a_{ij}) \in \mathbb{R}_{\max}^{n \times n}$, 
we define a weighted digraph $\mathcal{G}(A) = (N,E,w)$ as follows.
The sets of the nodes and arcs are $N=\{1,2,\dots,n\}$ and $E=\{ (i,j) \ |\ a_{ij} \neq \varepsilon \}$, respectively, and the weight function $w: E \to \mathbb{R}$ is defined by $w((i,j)) = a_{ij}$ for $(i,j) \in E$. A sequence of nodes $\mathcal{P} = (i_{0},i_{1},\dots,i_{\ell})$ is called a path if $(i_{k},i_{k+1}) \in E$ for $k=0,1,\dots,\ell-1$. It is called an $i_{0}$-$i_{\ell}$ path if its start and end should be specified. The sets of the nodes and arcs in $\mathcal{P}$ are denoted by $N(\mathcal{P})$ and $E(\mathcal{P})$, respectively. The number $\ell(\mathcal{P}) := \ell$ is called the length of $\mathcal{P}$. The sum $w(\mathcal{P}) := \sum_{k=0}^{\ell-1} w((i_{k},i_{k+1}))$ is called the weight of $\mathcal{P}$. A path $\mathcal{C} = (i_{0},i_{1},\dots,i_{\ell})$ with $i_{\ell} = i_{0}$ is called a circuit. In particular, if $i_{k} \neq i_{k'}$ for $1 \leq k < k' \leq \ell$, then $\mathcal{C}$ is called an elementary circuit. The length and weight of a circuit are defined similarly to a path. The mean weight of a circuit $\mathcal{C}$ is defined by $\mathrm{mw}(\mathcal{C}) := w(\mathcal{C}) / \ell(\mathcal{C})$. A union of disjoint elementary circuits is called a multi-circuit and its length (resp.~weight) is defined as the sum of the lengths (resp.~weights) of all circuits in it.
\par
For $A \in \mathbb{R}_{\max}^{n \times n}$ and a positive integer $k$, the $(i,j)$ entry of $A^{\otimes k}$ is identical to the maximum weight of all $i$-$j$ paths with lengths $k$ in $\mathcal{G}(A)$. We consider the formal matrix power series of the form
\begin{align*}
	A^* := I_{n} \oplus A \oplus A^{\otimes 2}\oplus \cdots.
\end{align*}
If there is no circuit with positive weight in $\mathcal{G}(A)$, then $A^{*}$ is computed as the finite sum 
\begin{align*}
	A^* = I_{n} \oplus A \oplus A^{\otimes 2}\oplus \cdots \oplus A^{\otimes n-1}.
\end{align*}
In this case, the $(i,j)$ entry of $A^{*}$ is the maximum weight of all $i$-$j$ paths.
\par
\subsection{Eigenvalues and eigenvectors}

For a matrix $A \in \mathbb{R}_{\max}^{n \times n}$, a scalar $\lambda$ is called an eigenvalue of $A$
if there exists a vector $\bm{x} \neq \mathcal{E}$ satisfying
\begin{align*}
	A \otimes \bm{x} = \lambda \otimes \bm{x}.
\end{align*}
This vector $\bm{x}$ is called an eigenvector of $A$ with respect to $\lambda$. Here, we summarize the results in the literature on the max-plus eigenvalue problem,
 e.g.,~\cite{Baccelli1992,Butkovic2010,Heidergott2006}.
\begin{proposition}\label{eigval1}
	For a matrix $A \in \mathbb{R}_{\max}^{n \times n}$, the maximum mean weight of all elementary circuits in $\mathcal{G}(A)$ is the maximum eigenvalue of $A$.
\end{proposition}
Let $\lambda(A)$ be the maximum eigenvalue of $A$.
A circuit in $\mathcal{G}(A)$ with mean weight $\lambda(A)$ is called critical. Nodes and arcs of a critical circuit is called critical nodes and critical arcs, respectively. The sets of all critical nodes and critical arcs are denoted by $N^{c}(A)$ and $E^{c}(A)$. The subgraph $\mathcal{G}^c(A) = (N^{c}(A), E^{c}(A))$ of $\mathcal{G}(A)$ is called the critical graph. 
\begin{proposition}\label{eigvec1}
	The $k$th column of $((-\lambda(A)) \otimes A)^*$ is an eigenvector of $A$ with respect to $\lambda(A)$ if and only if $k \in N^{c}(A)$.
\end{proposition}

A (univariate) polynomial in max-plus algebra has the form
\begin{align}
	f(t) = c_0 \oplus c_1 \otimes t \oplus c_2 \otimes t^{\otimes 2} \oplus \cdots \oplus c_n \otimes t^{\otimes n},   \label{eq:maxpoly}
\end{align}
where $c_0,c_1,\dots,c_n \in \mathbb{R}_{\max}$ and $c_{n} \neq \varepsilon$. It is a piecewise linear function on $\mathbb{R}_{\max}$.
Every polynomial can be factorized into a product of linear factors as follows:
\begin{align*}
	f(t) = c_{n} \otimes (t \oplus r_1)^{\otimes m_1} \otimes (t \oplus r_2)^{\otimes m_2} \otimes \cdots 
		\otimes (t \oplus r_p)^{\otimes m_p}.
\end{align*}
Then, $r_i$ and $m_i$ are referred to as the root of $f(t)$ and its multiplicity, respectively. At a root of $f(t)$, at least two terms on the right-hand side of \eqref{eq:maxpoly} attain the maximum value simultaneously. This property corresponds to the definition of tropical variety. In the graph of the piecewise linear function $f(t)$, roots are undifferentiable points of $f(t)$.
\par
The characteristic polynomial of $A \in \mathbb{R}_{\max}^{n \times n}$ is defined by
\begin{align*}
	\chi_A(t) := \det (A \oplus t \otimes I_{n}).
\end{align*}
As in conventional algebra, the characteristic polynomial of a matrix is closely related to the eigenvalue problem.
\begin{theorem}[\cite{Green1983}, Theorem 3]
For a matrix $A \in \mathbb{R}^{n \times n}_{\max}$, 
the maximum root of the characteristic polynomial is the maximum eigenvalue of $A$, which is also identical to the maximum mean circuit weight in $\mathcal{G}(A)$.
\end{theorem}
The roots of $\chi_{A}(t)$ are called algebraic eigenvalues of $A \in \mathbb{R}_{\max}^{n \times n}$~\cite{Akian2005}. Note that there are $n$ algebraic eigenvalues counting multiplicities.
\par
When we expand the polynomial $\chi_{A}(t)$, the coefficient of $t^{\otimes (n-k)}$ is the maximum weight of 
multi-circuits in $\mathcal{G}(A)$ with lengths $k$. For any fixed $\lambda \in \mathbb{R}$, a multi-circuit $\mathfrak{C}$ satisfying $\chi_{A}(\lambda) = w(\mathfrak{C}) \otimes \lambda^{\otimes (n-\ell(\mathfrak{C}))}$ is called a $\lambda$-maximal multi-circuit ($\lambda$-MMC). Let $\lambda_{1}, \lambda_{2}, \dots, \lambda_{p}$ be the finite roots of $\chi_{A}(t)$ in descending order. A sequence of multi-circuits $\mathfrak{C}_{0}, \mathfrak{C}_{1}, \dots, \mathfrak{C}_{p}$ is called the maximal multi-circuit sequence (MMCS) of $A$ if the following properties are satisfied.

\begin{enumerate}
    \item $\mathfrak{C}_{0} = \emptyset$.
    \item For $k=0,1,\dots,p$, $\mathfrak{C}_{k}$ is a $\lambda$-MMC for all $\lambda$ with $\lambda_{k+1} \leq \lambda \leq \lambda_{k}$. Here, $\lambda_{0}$ and $\lambda_{p+1}$ are considered to be $+\infty$ and $\varepsilon$, respectively.
\end{enumerate}
By definition, $\mathfrak{C}_{k}$ is the $\lambda_{k}$-MMC with the maximum length and $\lambda_{k+1}$-MMC with the minimum length.

Although finding all coefficients in the expansion of $\chi_{A}(t)$ is a difficult problem known as the best principal submatrix problem~\cite{Burkard2003-2}, Gassner and Klinz~\cite{Gassner2009} proposed a polynomial-time algorithm to obtain all roots of $\chi_{A}(t)$ using the parametric shortest path problem.

\begin{proposition}[\cite{Gassner2009}]
    For $A \in \mathbb{R}_{\max}^{n \times n}$ with $m$ finite entries, all roots of $\chi_{A}(t)$ and the MMCS of $A$ can be found in $O(n(m+n\log n))$ times of computation.
\end{proposition}

\section{CSR decomposition}

We consider a discrete event system defined by $A \in \mathbb{R}_{\max}^{n \times n}$ and $\bm{x}(0) \in \mathbb{R}_{\max}^{n}$ as follows:
\begin{align*}
    \bm{x}(t+1) = A \otimes \bm{x}(t), \quad t=0,1,\dots.
\end{align*}
Because $\bm{x}(t)$ can be written as $A^{\otimes t} \otimes \bm{x}(0)$, the power $A^{\otimes t}$ plays a significant role in deciding the asymptotic behavior of the system. If $A$ is irreducible, that is, the graph $\mathcal{G}(A)$ is strongly connected, then there exist integers $T$ and $\sigma$ such that 
\begin{align*}
    A^{\otimes (t+\sigma)} = \lambda(A)^{\otimes \sigma} \otimes A^{\otimes t}
\end{align*}
for all $t \geq T$. The minimum of such integer $T$ is denoted by $T(A)$ and called the transient of $A$. The transient $T(A)$ depends on the dimension $n$ of $A$, as well as on the values of the entries in $A$. Nachtigall~\cite{Nachtigall1997} introduced an expansion by at most $n$ terms:
\begin{align*}
    A^{\otimes t} = \bigoplus_{k=1}^{r} A_{k}^{\otimes t},
\end{align*}
where transients $T(A_{k})$ are $O(n^{2})$ for all $k$. Sergeev and Schneider~\cite{Sergeev2012} improved this expansion by introducing the CSR expansion of the form
\begin{align}
    A^{\otimes t} = \bigoplus_{k=1}^{r} \lambda_{k}^{\otimes t} \otimes C_{k} \otimes S_{k}^{\otimes t} \otimes R_{k}, 
    	\label{eq:csr1}
\end{align}
where $C_{k} \in \mathbb{R}_{\max}^{n \times n_{k}}$ and $ R_{k} \in \mathbb{R}_{\max}^{n_{k} \times n}$, and $S_{k} \in \{0,\varepsilon\}^{n_{k} \times n_{k}}$ is a permutation matrix, that is, each row and column has exactly one $0$ in its entries. Each $\lambda_{k}$ is called the growth rate of the corresponding term. An $O(n^{4} \log n)$ algorithm for CSR expansion is presented in~\cite{Sergeev2012}. Here, we briefly summarize the result. 

A max-plus matrix $A = (a_{ij}) \in \mathbb{R}_{\max}^{n \times n}$ is called visualized if $a_{ij} \leq \lambda(A)$ for all $i,j$. We can observe that $a_{ij} = \lambda(A)$ for all critical arcs $(i,j)$. Any matrix $A \in \mathbb{R}_{\max}^{n \times n}$ can be visualized into $A'$ by a diagonal matrix $D$ as $A' = D^{-1} \otimes A \otimes D$. It is easily verified that $\lambda(A) = \lambda(A')$ and $(A')^{\otimes t} = D^{-1} \otimes A^{\otimes t} \otimes D$. Therefore, we may assume that $A$ is a visualized matrix. In addition, it was shown in~\cite{Sergeev2009} that we may choose $D = \mathrm{diag}(\bm{d})$ such that $\bm{d}$ is any column of $(-\lambda(A) \otimes A)^{*}$. Specifically, $\bm{d} = (d_{1},d_{2},\dots,d_{n})^{\top}$ is the vector such that $d_{i}$ is the maximum weight of all paths from node $i$ to any fixed node $i_{0}$ in $\mathcal{G}(-\lambda(A) \otimes A)$.

In a strongly connected digraph $\mathcal{G}$, the cyclicity of $\mathcal{G}$ is the greatest common divisor of the lengths of all elementary circuits in $\mathcal{G}$. The cyclicity of a general graph is the least common multiple of the cyclicities of all strongly connected components. The cyclicity $\sigma(A)$ of a matrix $A \in \mathbb{R}_{\max}^{n \times n}$ is defined by the cyclicity of the critical graph $\mathcal{G}^{c}(A)$. The critical nodes are divided into equivalence classes with respect to the equivalence relation such that $i \sim j$ if and only if there exists an $i$-$j$ path whose length is a multiple of $\sigma(A)$. Let $\tilde{N}$ be the set of all equivalence classes. Let us consider the reduced graph $\tilde{\mathcal{G}} = (\tilde{N}, \tilde{E})$. The arcs are 
\begin{align*}
    \tilde{E} = \{ ([i],[j]) \mid \text{there exists an $i$-$j$ path whose length is $1 \bmod \sigma(A)$} \},
\end{align*}
where $[i]$ denotes the equivalence class containing node $i$. Then, $\tilde{\mathcal{G}}$ consists of disjoint elementary circuits.

Using the notations above, we describe the algorithm for the CSR expansion presented in~\cite{Sergeev2012}. Let $A \in \mathbb{R}_{\max}^{n \times n}$ be a visualized matrix with $\lambda(A) = 0$. The weights of all critical edges are $0$. We define $S \in \{0,\varepsilon\}^{|\tilde{N}| \times |\tilde{N}|}$ by the adjacent matrix of $\tilde{\mathcal{G}}$. In addition, we define $R \in \mathbb{R}_{\max}^{|\tilde{N}| \times n}$ by setting the $([i], j)$ entry to the maximum weight of all $i$-$j$ paths whose lengths are multiples of $\sigma(A)$. This is identical to the $(i,j)$ entry of $(A^{\otimes \sigma(A)})^{*}$. We similarly define $C \in \mathbb{R}_{\max}^{n \times |\tilde{N}|}$ by setting the $(i,[j])$ entry to the maximum weight of all $i$-$j$ paths whose lengths are multiples of $\sigma(A)$. Therefore, we obtain the first term of the CSR expansion of $A$ as $\lambda(A)^{\otimes t} \otimes C \otimes S^{\otimes t} \otimes R$. Subsequently, we delete all rows and columns of $A$ corresponding to $N^{c}(A)$, obtaining $A'$. Next, we obtain $\lambda(A')$ by the Karp Algorithm, visualize $(-\lambda(A')) \otimes A'$, and obtain another term $\lambda(A')^{\otimes t} \otimes C' \otimes (S')^{\otimes t} \otimes R'$ in the CSR expansion. We repeat this process until all rows and columns are deleted. Because the complexity of each step is at most $O(n^{3} \log n)$, the total computational complexity is at most $O(n^{4} \log n)$.

Several schemes for the CSR expansion were presented in~\cite{Merlet2014-2}. Nachtigall scheme comes from the original idea of the Nachtigall expansion~\cite{Nachtigall1997}. Hartmann-Arguelles scheme is based on the connectivity of the threshold graphs derived from the max-balanced scaling~\cite{Hartmann1999}. Another one is also based on the connectivity of threshold graphs, but they are determined by the mean weights of circuits.  The number $r$ of the terms in \eqref{eq:csr1}, as well as the minimum integer $T$ such that the equality \eqref{eq:csr1} holds for all $t \geq T$, differs depending on which scheme is used.

\section{Our algorithm}

In this study, we propose a new scheme for the CSR expansion of $A \in \mathbb{R}_{\max}^{n \times n}$, deriving an $O(n(m+n \log n))$ algorithm, where $m$ is the number of finite entries in $A$. In particular, the complexity is $O(n^{3})$ for dense matrices. This algorithm seems to be fast enough for this problem because the asymptotic growth rate $\lambda(A)$ is the maximum mean weight of all circuits in $\mathcal{G}(A)$, which is found by Karp Algorithm~\cite{Karp1978} in $O(mn)$ times of computation. Our contribution is summarized into the following three points.

.

\begin{enumerate}
    \item To obtain all values $\lambda_{k}$ in the CSR expansion, we use the characteristic polynomial $\chi_{A}(t)$ of the matrix $A$. In particular, we show that every $\lambda_{k}$ is a root of $\chi_{A}(t)$. We apply the algorithm of Gassner and Klinz~\cite{Gassner2009} to obtain growth rates $\lambda_{1}, \lambda_{2}, \dots, \lambda_{r}$ in the CSR expansion and matrices $A_{1}, A_{2}, \dots, A_{r}$ such that $A_{k+1}$ is a submatrix of $A_{k}$ with $\lambda(A_{k+1}) = \lambda_{k+1}$. The computational complexity is $O(n(m+n \log n))$.
    \item We visualize all matrices $A_{1}, A_{2}, \dots, A_{r}$ in $O(n(m+n \log n))$ times of computation in total by applying the Dijkstra Algorithm successively.
    \item Let $\ell$ be the length of the critical circuit $\mathcal{C}$ in $\mathcal{G}(A)$. We propose a graph extension to compute the critical rows and columns of $(A^{\otimes \ell})^{*}$ in $O(\ell(m+n\log n))$ times of computation. 
\end{enumerate}

In the following subsections, we discuss each point.

\subsection{Decomposition based on roots of characteristic polynomials} \label{sec:41}

Let $\lambda_{1}, \lambda_{2}, \dots, \lambda_{p}$ be the finite roots of the characteristic polynomial of $A \in \mathbb{R}_{\max}^{n \times n}$ in descending order. Recall that we can obtain the MMCS $\mathfrak{C}_{0}, \mathfrak{C}_{1}, \dots, \mathfrak{C}_{p}$ simultaneously. We define a partition of the node set $N$ by Algorithm \ref{alg:partition}. In each iteration, we take a circuit $\mathcal{C}$ from $\mathfrak{C}_{k}$ and check whether $\mathcal{C}$ intersects another circuit that was previously managed. If it does, the nodes of $\mathcal{C}$ are added to the last group obtained so far. Otherwise, $N(\mathcal{C})$ itself creates the new group.
Because $p \leq n$ and the inner for-loop takes $O(|N(\mathcal{C})|)$ times of computation, the complexity of the algorithm is $O(n(m+n \log n))$. 

 \begin{algorithm}
 \caption{Node partition}
 \label{alg:partition}
 \begin{algorithmic}
 \AlgInput{The MMCS $\mathfrak{C}_{0}, \mathfrak{C}_{1}, \dots, \mathfrak{C}_{p}$ corresponding to $A \in \mathbb{R}_{\max}^{n \times n}$}
 \AlgOutput{A partition of nodes $N_{1}, N_{2}, \dots, N_{r}$}
 \State{Define $N_{0} := \emptyset, U_{0} := \emptyset$ and $r := 0$}
 \For{$k = 1,2,\dots,p$}
 \For{$\mathcal{C} \in \mathfrak{C}_{k}$}
 \If{$N(\mathcal{C}) \cap U_{r} \neq \emptyset$}
 \State{Update $N_{r} := N_{r} \cup (N(\mathcal{C}) \setminus U_{r})$ and $U_{r} := U_{r} \cup N(\mathcal{C})$}
 \Else
 \State{Set $N_{r+1} := N(\mathcal{C})$ and $U_{r+1} := U_{r} \cup N(\mathcal{C})$}
 \State{Update $r := r+1$}
 \EndIf
 \EndFor
 \EndFor
 \State{Update $N_{r} := N_{r} \cup (N \setminus U_{r})$ and $U_{r} := N$}\\
 \Return $N_{1}, N_{2}, \dots, N_{r}$
 \end{algorithmic}
 \end{algorithm}

If a partition $N_{s}$ is initially created by a circuit $\mathcal{C}$ in $\mathfrak{C}_{k}$, then we denote such integer $k$ by $k(s)$. This circuit $\mathcal{C}$ is called quasi-critical. We prove that the CSR expansion of $A$ is expressed by $r$ terms whose growth rates are $\lambda_{k(1)}, \lambda_{k(2)}, \dots, \lambda_{k(r)}$. For $s=1,2,\dots,r$, we can observe that $U_{s} = \bigcup_{s' \leq s} N_{s'}$ by definition. Further, let $V_{s} = N \setminus U_{s-1}$. We first present several lemmas for $\lambda$-MMC.

\begin{lemma} \label{lem:mcc1}
    For $k=1,2,\dots,p$, the mean weight of any circuit in $\mathfrak{C}_{k}$ is at least $\lambda_{k}$.
\end{lemma}

\proof
    Let $\mathcal{C}$ be a circuit in $\mathfrak{C}_{k}$. Then, $\mathfrak{C}' = \mathfrak{C}_{k} \setminus \{\mathcal{C}\}$ is a multi-circuit. By the $\lambda_{k}$-maximality of $\mathfrak{C}_{k}$, we have
    \begin{align*}
        w(\mathfrak{C}_{k}) + \lambda_{k} (n-\ell(\mathfrak{C}_{k}))
    &\geq w(\mathfrak{C}') + \lambda_{k} (n-\ell(\mathfrak{C}')) \\ 
    &= w(\mathfrak{C}_{k}) + \lambda_{k} (n-\ell(\mathfrak{C}_{k})) - (w(\mathcal{C}) - \lambda_{k} \ell(\mathcal{C})).
    \end{align*}
    This implies $\mathrm{mw}(\mathcal{C}) \geq \lambda_{k}$.
\endproof

\begin{lemma} \label{lem:mcc2}
    For $s=1,2,\dots,r$, let $\mathcal{C}$ be the quasi-critical circuit in $N_{s}$. Then $\mathrm{mw}(\mathcal{C}) = \lambda_{k(s)}$. Furthermore, in the subgraph $\mathcal{G}(V_{s})$ induced by the node set $V_{s}$, the maximum mean weight of circuits is $\lambda_{k(s)}$.
\end{lemma}

\proof
    By definition, $\mathcal{C}$ is a circuit in $\mathfrak{C}_{k(s)}$. From Lemma \ref{lem:mcc1}, we have $\mathrm{mw}(\mathcal{C}) \geq \lambda_{k(s)}$. 
    On the other hand, let $\mathcal{C}'$ be any circuit in $V_{s}$. Because all nodes of $\mathfrak{C}_{k(s)-1}$ are contained in $U_{s-1}$ by Algorithm \ref{alg:partition}, $\mathfrak{C}' = \mathfrak{C}_{k(s)-1} \cup \{\mathcal{C}'\}$ is also a multi-circuit. By the $\lambda_{k(s)}$-maximality of $\mathfrak{C}_{k(s)-1}$, we have 
    \begin{align*}
        &w(\mathfrak{C}_{k(s)-1}) + \lambda_{k(s)} (n-\ell(\mathfrak{C}_{k(s)-1})) \\
    \geq&\, w(\mathfrak{C}') + \lambda_{k(s)} (n-\ell(\mathfrak{C}')) \\ 
    =&\, w(\mathfrak{C}_{k(s)-1}) + \lambda_{k(s)} (n-\ell(\mathfrak{C}_{k(s)-1})) + (w(\mathcal{C}') - \lambda_{k} \ell(\mathcal{C}')).
    \end{align*}
    This implies $\mathrm{mw}(\mathcal{C}') \leq \lambda_{k(s)}$. Thus, $\mathrm{mw}(\mathcal{C}) = \lambda_{k(s)}$, and it is the maximum mean weight of circuits in $\mathcal{G}(V_{s})$.
\endproof

\begin{lemma} \label{lem:mcc3}
    For $k=1,2,\dots,p$, if a circuit $\mathcal{C}$ and $\mathfrak{C}_{k}$ have a common node, then $\mathcal{C}$ intersects some circuit $\mathcal{C}'$ in one of $\mathfrak{C}_{1}, \mathfrak{C}_{2}, \dots, \mathfrak{C}_{k}$ such that $\mathrm{mw}(\mathcal{C}') \geq \mathrm{mw}(\mathcal{C})$.
\end{lemma}

\proof
    The case $k=1$ is trivial because $\mathfrak{C}_{1}$ consists of critical circuits. For $k \geq 2$, let $k' \leq k$ be the maximum integer such that $\mathcal{C}$ and $\mathfrak{C}_{k'-1}$ have no common node. By the $\lambda_{k'}$-maximality of $\mathfrak{C}_{k'-1}$, we have
    \begin{align*}
        &w(\mathfrak{C}_{k'-1}) + \lambda_{k'} (n-\ell(\mathfrak{C}_{k'-1})) \\
    \geq&\, w(\mathfrak{C}_{k'-1} \cup \{\mathcal{C}\}) + \lambda_{k'} (n-\ell(\mathfrak{C}_{k'-1} \cup \{\mathcal{C}\})) \\ 
    =&\, w(\mathfrak{C}_{k'-1}) + \lambda_{k'} (n-\ell(\mathfrak{C}_{k'-1})) + (w(\mathcal{C}) - \lambda_{k'} \ell(\mathcal{C})).
    \end{align*}
    This implies that $\mathrm{mw}(\mathcal{C}) \leq \lambda_{k'}$. By the choice of $k'$, there exists a circuit $\mathcal{C}' \in \mathfrak{C}_{k'}$ that intersects $\mathcal{C}$. By Lemma \ref{lem:mcc1}, $\mathrm{mw}(\mathcal{C}') \geq \lambda_{k'}$. Hence, we have $\mathrm{mw}(\mathcal{C}') \geq \mathrm{mw}(\mathcal{C})$.
\endproof

\begin{lemma} \label{lem:int}
    Let $\{a_{1}, a_{2}, \dots, a_{\ell}\}$ be a set of integers. There exists a subset $I \subset \{1,2,\dots,\ell\}$ such that $\sum_{i \in I} a_{i}$ is a multiple of $\ell$.
\end{lemma}

\proof
    For $j = 1,2,\dots, \ell$, let $b_{j} = \sum_{i=1}^{j} a_{i}$. If some $b_{j}$ is a multiple of $\ell$, the assertion of the lemma is trivial. Otherwise, by the pigeonhole principle, there are two integers $j_{1}$ and $j_{2}$, where $j_{1} < j_{2}$, such that $b_{j_{1}} = b_{j_{2}} \bmod \ell$. Therefore, $b_{j_{2}} - b_{j_{1}} = \sum_{i=j_{1}+1}^{j_{2}} a_{i}$ is a multiple of $\ell$.
\endproof

Now, we state that any path in $\mathcal{G}(A)$ with the maximum weight contains a quasi-critical circuit if the length is more than $(n+1)n$.

\begin{proposition} \label{prop:longpath}
    Let $A \in \mathbb{R}_{\max}^{n \times n}$. For any $t > (n+1)n$ and nodes $i$ and $j$, let $\Pi_{i,j,s}^{t}$ be a set of $i$-$j$ paths with lengths $t$ that pass through some nodes in $U_{s}$. Then, a path $\mathcal{P} \in \Pi_{i,j,s}^{t}$ with maximum weight contains a quasi-critical circuit in $U_{s}$ as its subpath. In particular, if $\mathcal{P}$ does not pass through $U_{s-1}$, then it contains a quasi-critical circuit in $N_{s}$.
\end{proposition}

\proof
    Let $\mathcal{P}$ be any path in $\Pi_{i,j,s}^{t}$. We show that there exists a path $\mathcal{P}' \in \Pi_{i,j,s}^{t}$ such that $\mathcal{P}'$ contains a quasi-critical circuit in $U_{s}$ and $w(\mathcal{P}') \geq w(\mathcal{P})$.

    If $\mathcal{P}$ contains a quasi-critical circuit in $U_{s}$, the assertion is trivial. Otherwise, let $\mathcal{P}_{1} = \mathcal{P}$ and $\mathcal{C}_{1}$ be a circuit with the maximum mean weight among those that intersect, but are not contained in, $\mathcal{P}_{1}$ and are contained in $\mathfrak{F}_{s} = \mathfrak{C}_{1} \cup \mathfrak{C}_{2} \cup \dots \cup \mathfrak{C}_{k(s)}$. If there are two or more circuits, we take one from $\mathfrak{C}_{k}$ with the smallest $k$. Because $\mathcal{P}$ passes through $U_{s}$, such a circuit $\mathcal{C}_{1}$ always exists. Let $i_{1}$ be a node where $\mathcal{P}$ intersects $\mathcal{C}_{1}$ firstly. The elementary $i$-$i_{1}$ and $i_{1}$-$j$ paths contained in $\mathcal{P}_{1}$ are denoted by $\mathcal{P}^{-}_{1}$ and $\mathcal{P}^{+}_{1}$, respectively. The other part of $\mathcal{P}_{1}$, denoted by $\mathfrak{P}_{1}$, consists of circuits, where identical circuits are counted many times. We have
    \begin{align*}
        \ell(\mathcal{P}^{-}_{1}) + \ell(\mathcal{P}^{+}_{1}) + \ell(\mathcal{C}_{1}) \leq 2n.
    \end{align*}
    The total length of the circuits in $\mathfrak{P}_{1}$ is more than $(n-1)n$. Thus, $\mathfrak{P}_{1}$ contains at least $n$ circuits. By Lemma \ref{lem:int}, we can choose a collection of circuits in $\mathfrak{P}_{1}$ to ensure that the sum of the lengths is a multiple of $\ell(\mathcal{C}_{1})$, say $m_{1} \ell(\mathcal{C}_{1})$. Subsequently, we remove these circuits from $\mathcal{P}_{1}$ and append $m_{1}$ times of $\mathcal{C}_{1}$, obtaining a new $i$-$j$ path $\mathcal{P}_{2}$ with length $t$. We can observe that $\mathcal{P}_{2}$ also passes through $U_{s}$ because $\mathcal{C}_{1} \in \mathfrak{F}_{s}$. In addition, we can verify that $w(\mathcal{P}_{2}) \geq w(\mathcal{P}_{1})$. For all circuits $\mathcal{C}' \in \mathfrak{P}_{1}$, we have $\mathrm{mw}(\mathcal{C}') \leq \mathrm{mw}(\mathcal{C}_{1})$. Indeed, if $V(\mathcal{C}') \cap U_{s} \neq \emptyset$ , then there exists $\mathcal{C}'' \in \mathfrak{F}_{s}$ such that $\mathrm{mw}(\mathcal{C}') \leq \mathrm{mw}(\mathcal{C}'')$ by Lemma \ref{lem:mcc3}. By the maximality in choosing $\mathcal{C}_{1}$, we can observe that $\mathrm{mw}(\mathcal{C}'') \leq \mathrm{mw}(\mathcal{C}_{1})$, yielding $\mathrm{mw}(\mathcal{C}') \leq \mathrm{mw}(\mathcal{C}_{1})$. On the other hand, if $V(\mathcal{C}') \cap U_{s} = \emptyset$, then $\mathrm{mw}(\mathcal{C}') \leq \lambda_{k(s+1)} \leq \lambda_{k(s)} \leq \mathrm{mw}(\mathcal{C}_{1})$ according to Lemma \ref{lem:mcc1} and \ref{lem:mcc2}. Thus, exchange of circuits in $\mathcal{P}_{1}$ in the above manner does not decrease the weight of the path. 

    If $\mathcal{C}_{1}$ is not a quasi-critical circuit in $U_{s}$, we next obtain a circuit $\mathcal{C}_{2}$ with the maximum mean weight among those that intersect $\mathcal{C}_{1}$. By the maximality in choosing $\mathcal{C}_{1}$, we can observe that $\mathcal{C}_{2}$ does not intersect $\mathcal{P}$. We similarly define $\mathcal{P}^{-}_{2}$, $\mathcal{P}^{+}_{2}$ and $\mathfrak{P}_{2}$. We obtain
    \begin{align*}
        \ell(\mathcal{P}^{-}_{2}) + \ell(\mathcal{P}^{+}_{2}) + \ell(\mathcal{C}_{2}) \leq 2n.
    \end{align*}
    Thus, we can choose a collection of circuits to ensure that the sum of the lengths is $m_{2} \ell(\mathcal{C}_{2})$. We remove these circuits from $\mathcal{P}_{2}$ and append $m_{2}$ times $\mathcal{C}_{2}$, obtaining a new $i$-$j$ path $\mathcal{P}_{3}$ with length $t$ and weight $w(\mathcal{P}_{3}) \geq w(\mathcal{P}_{2})$. Repeating this process, we finally reach a quasi-critical circuit $\mathcal{C}_{q}$ in $U_{s}$ and obtain a desired path $\mathcal{P}' = \mathcal{P}_{q}$.
\endproof

For $s=1,2,\dots,r$, let $A(V_{s})$ be the matrix whose rows and columns are restricted to $V_{s}$. Proposition \ref{prop:longpath} suggests that we can consider the expansion
\begin{align}
    A^{\otimes t} = \bigoplus_{s=1}^{r} W^{t}_{s},  \label{eq:nach}
\end{align}
where $(i,j)$ entry of $W^{t}_{s}$ is the maximum weight of $i$-$j$ paths in $\mathcal{G}(A(V_{s}))$ whose length is $t$ and that contains a quasi-critical circuit $\mathcal{C}$ in $N_{s}$.

\subsection{Visualization of matrices}

For $s=1,2,\dots,r$, let $n_{s}$ and $m_{s}$ be the number of nodes and arcs in $\mathcal{G}(A(V_{s}))$. Each matrix $A_{s}$ can be visualized in the complexity $O(m_{s}n_{s})$ by the Bellman-Ford Algorithm to obtain the shortest (maximum weight) path tree with possibly positive arc weights. Here, larger weights of arcs are considered to be shorter by reversing the signs of the weights of the arcs. Therefore, the total complexity to visualize $A_{s}$ for all $s=1,2,\dots,r$ is at most $O(mn^{2})$ in total. However, in Algorithm \ref{alg:visualize}, we can efficiently visualize them in backward order, i.e., from $A_{r}$ to $A_{1}$, by using the visualization of smaller matrices. We remark that the Dijkstra Algorithm can be slightly extended.

\begin{lemma} \label{lem:dijkstra}
    Let $\mathcal{G} = (N,E)$ be a digraph with length function $l:E \to \mathbb{R}$, and suppose that there is no negative circuit in $\mathcal{G}$ with respect to $l$. Let $s$ be a node of $\mathcal{G}$. If only the arcs incident to $s$ have negative lengths, then the shortest path tree from the root $s$ can be obtained in $O(|E|+|N|\log |N|)$ times of computation by the Dijkstra Algorithm.
\end{lemma}

\proof
    Because the arcs with negative lengths appear only in the first iteration, the minimum value of the labels on the nodes remained unchanged throughout the algorithm.
\endproof

\begin{algorithm}
 \caption{Visualization}
 \label{alg:visualize}
 \begin{algorithmic}
 \AlgInput{The partition of nodes $N_{1}, N_{2}, \dots, N_{r}$ obtained by Algorithm \ref{alg:partition}}
 \AlgOutput{Visualization of matrices $A(V_{1}), A(V_{2}), \dots, A(V_{r})$ and the vectors $\bm{d}_{1}, \bm{d}_{2}, \dots, \bm{d}_{r}$ inducing them}
 \State{Define $N' := \emptyset$ and $d_{j} := 0$ for all $j \in N$}
 \For{$s = r,r-1,\dots,1$}
 \State{Update $b((j_{1},j_{2})) := b((j_{1},j_{2})) + (\lambda_{k(s+1)} - \lambda_{k(s)})$ for $(j_{1},j_{2}) \in (N' \times N') \cap E$}
 \For{$i \in N_{s}$}
 \State{Set $b((i,j)) := (-\lambda_{k(s)}) + [A]_{ij} + d_{j}$ for $j \in N'$ with $(i,j) \in E$}
 \State{Set $b((j,i)) := -d_{j} + (-\lambda_{k(s)}) + [A]_{ij}$ for $j \in N'$ with $(j,i) \in E$}
 \State{Set $b((i,i)) := (-\lambda_{k(s)}) + [A]_{ii}$}
 \State{Update $N' := N' \cup \{i\}$}
 \State{In graph $\mathcal{G}(b(N'))$, find the maximum weight path tree going into the root $i$}
 \State{Set $R \subset N'$ be the nodes that are reachable to $i$}
 \State{Set $w_{j}$ to the maximum weight of $j$-$i$ path for $j \in R$}
 \State{Update $b((j_{1},j_{2})) := -w_{j_{1}} + b((j_{1},j_{2})) + w_{j_{2}}$ for $(j_{1},j_{2}) \in (R \times R) \cap E$}
 \State{Set $w^{*} = -\max(0, \max\{ -w_{j_{1}} + b((j_{1},j_{2})) \mid (j_{1},j_{2}) \in (R \times (N' \setminus R)) \cap E\})$}
 \State{Set $w_{j} = w^{*}$ for $j \in N' \setminus R$}
 \State{Update $b((j_{1},j_{2})) := -w_{j_{1}} + b((j_{1},j_{2})) + w_{j_{2}}$ for $(j_{1},j_{2}) \in (R \times (N' \setminus R)) \cap E$}
 \State{Update $d_{j} := d_{j} + w_{j}$ for $j \in R$}
 \EndFor
 \State{Set $[A'_{s}]_{j_{1}j_{2}} := b((j_{1},j_{2}))$ for $(j_{1},j_{2}) \in (N' \times N') \cap E$}
 \State{Set $[\bm{d}_{s}]_{j} := d_{j}$ for $j \in N'$}
 \EndFor\\
 \Return $A'_{r}, A'_{r-1}, \dots, A'_{1}$ and $\bm{d}_{r}, \bm{d}_{r-1}, \dots, \bm{d}_{1}$
 \end{algorithmic}
 \end{algorithm}

Consider Algorithm \ref{alg:visualize}. In the iteration for $i \in N_{s}$, we first compute the weights of the arcs incident to the new node $i$. This reflects the scaling by $d_{j}$ cumulated in the previous iterations. We obtain the maximum weight path tree of the graph $\mathcal{G}(b(N'))$, where $b(N')$ is the $|N'| \times |N'|$ matrix with entries $b((j_{1},j_{2}))$. Note that we consider the single sink problem with sink $i$, which is identical to the single source problem by reversing all arcs. By updating $b((j_{1},j_{2})) := -w_{j_{1}} + b((j_{1},j_{2})) + w_{j_{2}}$, all values $b((j_{1},j_{2}))$ become non-positive because $w_{j}$ is a feasible potential. Later, we manage $b((j_{1},j_{2}))$ for $j_{2} \not\in R$. Because $b((j_{1},j_{2})) = \varepsilon$ for $(j_{1},j_{2}) \in (N' \setminus R) \times R$, we only need to update $b((j_{1},j_{2}))$ for $(j_{1},j_{2}) \in R \times (N' \setminus R)$ to ensure that $b((j_{1},j_{2})) \leq 0$ for all $j_{1},j_{2} \in N'$. When we move from the phase $s+1$ to $s$, we update by $b((j_{1},j_{2})) := b((j_{1},j_{2})) + (\lambda_{k(s+1)} - \lambda_{k(s)})$ to manage the matrix $(-\lambda_{s}) \otimes A_{s}$. Because $\lambda_{k(s+1)} \leq \lambda_{k(s)}$, this update does not increase the value $b((j_{1},j_{2}))$. Therefore, in each iteration, the weights of the arcs in $\mathcal{G}(b(N'))$ are non-positive, except those incident to the new node $i$. By Lemma \ref{lem:dijkstra}, we can apply the Dijkstra Algorithm in finding the maximum weight path tree. Because each node $i$ is appended to $N'$ exactly once, the total computational complexity of Algorithm \ref{alg:visualize} is $O(n(m+n\log n))$.

\subsection{Computing the CSR decomposition}

In the expansion \eqref{eq:nach}, we will derive the decomposition $W^{t}_{s} = \lambda_{k(s)}^{\otimes t}  \otimes C_{s} \otimes S_{s}^{\otimes t} \otimes R_{s}$ for $s=1,2,\dots,r$. Because this is considered for each submatrix $A(V_{s})$ separately, we may assume that $A \in \mathbb{R}_{\max}^{n \times n}$ is a visualized matrix with $\lambda(A) = 0$. Furthermore, without loss of generality, let $\mathcal{C}^{*} = \{1,2,\dots,\ell\}$ be the quasi-critical circuit in $N_{1}$, which is a critical circuit in $\mathcal{G}(A)$. We denote $W^{t}_{1}$ by $W^{t}$.

\begin{proposition} \label{prop:pathdiv}
    Let $d_{ij}$ be the maximum weight of all $i$-$j$ paths in $\mathcal{G}(A)$ whose lengths are multiples of $\ell$. For $t > (n+1)n$, we have 
    \begin{align*}
        [W^{t}]_{ij} \leq d_{ik} + d_{k'j},
    \end{align*}
    for some $k$ and $k'$ such that $k' = t+k \bmod \ell$.
\end{proposition}

\proof
Let $\mathcal{P} \in \Pi_{i,j,1}^{t}$ be the path attaining the maximum weight $[W^{t}]_{ij}$. By Proposition \ref{prop:longpath}, $\mathcal{P}$ contains $\mathcal{C}^{*}$. Subsequently, $\mathcal{P}$ can be divided into three parts: the $i$-$k$ path $\mathcal{P}_{1}$ whose length is a multiple of $\ell$, the $k$-$k'$ path along $C^{*}$, and the $k'$-$j$ path $\mathcal{P}_{3}$ whose length is a multiple of $\ell$, where $k,k' \in \{1,2,\dots,\ell\}$. We obtain $w(\mathcal{P}_{1}) \leq d_{ik}$ and $w(\mathcal{P}_{3}) \leq d_{k'j}$. The weights of all arcs in $\mathcal{C}^{*}$ are $0$ because $A$ is visualized. Hence, 
\begin{align*}
    [W^{t}]_{ij} = w(\mathcal{P}) = w(\mathcal{P}_{1}) + w(\mathcal{P}_{3}) \leq d_{ik} + d_{k'j}.
\end{align*}    
Furthermore, because $\ell(\mathcal{P}) = \ell(\mathcal{C}^{*}) \bmod \ell$, we have $t = k'-k \bmod \ell$. Thus, we have proved the proposition.
\endproof

\begin{proposition} \label{prop:CSR}
    Let $C \in \mathbb{R}_{\max}^{n \times \ell}$ and $R \in \mathbb{R}_{\max}^{\ell \times n}$ be the matrices whose $(i,j)$ entries are $d_{ij}$. Furthermore, let $S \in \{0,\varepsilon\}^{\ell \times \ell}$ whose $(i,j)$ entry is $0$ if and only if $(i,j) \in E(\mathcal{C}^{*})$. Then, $S$ is periodic with $S^{\otimes \ell} = I_{\ell}$ and
    \begin{align}
        W^{t} = C \otimes S^{\otimes t} \otimes R \label{eq:prf43}
    \end{align}
    for $t \geq 2n^{2}$.
\end{proposition}

\proof
By proposition \ref{prop:pathdiv}, we can observe that 
\begin{align*}
    [W^{t}]_{ij} \leq [C]_{ik} \otimes [S^{\otimes t}]_{kk'} \otimes [R]_{k'j} \leq [C \otimes S^{\otimes t} \otimes R]_{ij} 
\end{align*}
for all $i,j \in N$. To prove the opposite inequality, for any $i,j \in N$ and $k,k' \in \{1,2,\dots,\ell\}$ with $k' = t+k \bmod \ell$, let $\mathcal{P}_{1}$ and $\mathcal{P}_{2}$ be the maximum weight $i$-$k$ path and $k'$-$j$ path whose lengths are multiple of $\ell$, respectively. Because these correspond to the maximum weight paths in $\mathcal{G}(A^{\otimes \ell})$, the lengths of both $\mathcal{P}_{1}$ and $\mathcal{P}_{2}$ are at most $\ell n \leq n^{2}$. Thus, if $t \geq 2n^{2}$, we can construct an $i$-$j$ path $\mathcal{P}$ with length $t$ by a composition of $\mathcal{P}_{1}$, the $k$-$k'$ path with length $t - \ell(\mathcal{P}_{1}) - \ell(\mathcal{P}_{2})$ along $\mathcal{C}^{*}$, and $\mathcal{P}_{2}$. Therefore, $w(\mathcal{P}) \leq [W^{t}]_{ij}$. Because all terms appearing in the $(i,j)$ entry of the right-hand side of \eqref{eq:prf43} can similarly be obtained as $i$-$j$ paths, we can conclude that
\begin{align*}
    [W^{t}]_{ij} \geq [C \otimes S^{\otimes t} \otimes R]_{ij}. 
\end{align*}
\endproof

\begin{corollary}[CSR expansion]		\label{coro:csr}
    For $t \geq 2n^{2}$, $A\in \mathbb{R}_{\max}^{n \times n}$ can be expanded as 
    \begin{align*}
        A^{\otimes t} = \bigoplus_{s=1}^{r} \lambda_{k(s)}^{\otimes t} \otimes C_{s} \otimes S_{s}^{\otimes t} \otimes R_{s}.
    \end{align*}
    Here, $S_{s}$ is a $(0,\varepsilon)$-matrix whose $(i,j)$ entry is $0$ if $(i,j) \in E(\mathcal{C}_{s})$ for the quasi-critical circuit $\mathcal{C}_{s}$ in $N_{s}$.
\end{corollary}

\begin{remark}
The CSR expansion can be improved in terms of the cyclicity theorem~\cite{Cohen1985}. As described in Section 3, the decomposition \eqref{eq:prf43} can be expressed as $W^{t} = \tilde{C} \otimes \tilde{S}^{\otimes t} \otimes \tilde{R}$, where $\tilde{C}, \tilde{S}$ and $\tilde{R}$ are defined in terms of the reduced graph $\tilde{\mathcal{G}}$. In particular, the period of $\tilde{S}$ is the cyclicity of $A$. Thus, we obtain the CSR expansion
    \begin{align*}
        A^{\otimes t} = \bigoplus_{k} \lambda_{k}^{\otimes t} \otimes \tilde{C}_{k} \otimes \tilde{S}_{k}^{\otimes t} \otimes \tilde{R}_{k}.
    \end{align*}
    Here, $k$ takes over the distinct integers among $k(1),k(2),\dots,k(r)$. For each $k$, the decomposition $\tilde{C}_{k} \otimes \tilde{S}_{k}^{\otimes t} \otimes \tilde{R}_{k}$ is based on the reduced graph of $\mathcal{G}((-\lambda_{k}) \otimes A(V_{s}))$, where $s$ is the minimum integer such that $k=k(s)$. 
\end{remark}

We describe a way to obtain $C_{s}$ and $R_{s}$ efficiently. We consider $A \in \mathbb{R}_{\max}^{n \times n}$ as in the beginning of this subsection. From $\mathcal{G}(A) = (N,E,w)$, we construct the extended graph $\hat{\mathcal{G}} = (\hat{N},\hat{E},\hat{w})$ as follows. The nodes of $\hat{\mathcal{G}}$ are expressed by the set $\hat{N} = \{ (i,k) \mid i \in N, k = \{0,1,\dots,\ell-1\} \}$. The arc set is defined by 
\begin{align*}
    \hat{E} = \{ ((i,k),(j,k+1) \mid (i,j) \in E, k=0,1,\dots,\ell-2 \} \cup \{ ((i,\ell-1),(j,0) \mid (i,j) \in E \}.
\end{align*}
The weights of the arcs in $\hat{E}$ inherit those in $E$. Because $A$ is visualized, all arcs in $\hat{E}$ have non-positive weight. Thus, we apply the Dijkstra Algorithm to obtain the maximum weight path tree whose root is $(1,0) \in \hat{N}$. The computational complexity is 
\begin{align*}
    O(|\hat{E}| + |\hat{N}| \log |\hat{N}|) = O(\ell |E| + \ell |N| \log \ell |N|) = O(\ell (|E| + |N| \log |N| )),
\end{align*}
The maximum weight of $(1,0)$-$(j,k)$ path in $\hat{\mathcal{G}}$ is identical to the maximum weight of $1$-$j$ path in $\mathcal{G}(A)$ whose length is $k \bmod \ell$. Furthermore, it is identical to the maximum weight of $(k+1)$-$j$ path whose length is a multiple of $\ell$ since the weights of edges in $\mathcal{C}^*$ are 0. Thus, we have obtained the matrix $R$ in Proposition \ref{prop:CSR}. The matrix $C$ can be obtained in a similar way.

Suppose that the visualized matrices $A'_{1}, A'_{2}, \dots, A'_{r}$ are obtained from matrices $A(V_{1}), A(V_{2}), \dots, A(V_{r})$, respectively. For $s=1,2,\dots,r$, let $E_{s}$ be the arc set of $\mathcal{G}(A'_{s})$, while $V_{s}$ is the node set, and $\ell_{s}$ be the length of the quasi-critical circuit in $N_{s}$. We can observe that $\sum_{s=1}^{m} \ell_{s} \leq n$. From the above observation, we can obtain $C_{1}, C_{2}, \dots, C_{r}$ and $R_{1}, R_{2}, \dots, R_{r}$ in at most
\begin{align*}
    \sum_{s=1}^{r} O(\ell_{s} (|E_{s}| + |V_{s}| \log |V_{s}|)) 
    \leq O(m + n\log n) \cdot \sum_{s=1}^{r} O(\ell_{s})
    \leq O(n(m+n \log n))
\end{align*}
times of computation. Summarizing the results in this section, the computational complexity for the CSR expansion of $A \in \mathbb{R}_{\max}^{n \times n}$ with $m$ finite entries is at most $O(n(m+n \log n))$. 

\subsection{Example}

\begin{figure}[tbp]
 \centering
 \label{fig:exm2}\includegraphics[width=0.5\textwidth]{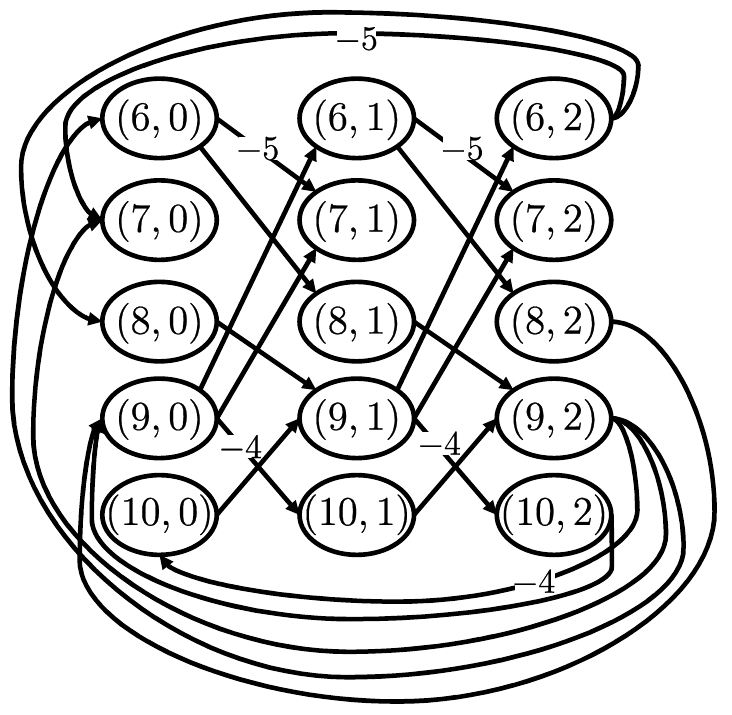}
 \caption{Extended graph $\hat{G}$. The weights of the arcs are $0$ if they are not specified.}
 \end{figure}

Let us consider a matrix
\begin{align*}
    A = \begin{pmatrix}
        \varepsilon & 7 & \varepsilon & \varepsilon & \varepsilon & \varepsilon & \varepsilon & \varepsilon & \varepsilon & \varepsilon \\
        9 & \varepsilon & 8 & 3 & 7 & \varepsilon & \varepsilon & \varepsilon & \varepsilon & \varepsilon \\
        8 & \varepsilon & \varepsilon & \varepsilon & \varepsilon & \varepsilon & \varepsilon & \varepsilon & \varepsilon & \varepsilon \\
        \varepsilon & \varepsilon & \varepsilon & 6 & 2 & 5 & \varepsilon & \varepsilon & \varepsilon & \varepsilon \\
        \varepsilon & \varepsilon & 5 & \varepsilon & \varepsilon & \varepsilon & \varepsilon & \varepsilon & \varepsilon & \varepsilon \\
        \varepsilon & \varepsilon & \varepsilon & \varepsilon & \varepsilon & \varepsilon & 1 & 6 & \varepsilon & \varepsilon \\        
        \varepsilon & \varepsilon & \varepsilon & \varepsilon & 2 & \varepsilon & \varepsilon & \varepsilon & \varepsilon & \varepsilon \\
        \varepsilon & \varepsilon & \varepsilon & \varepsilon & \varepsilon & \varepsilon & \varepsilon & \varepsilon & 2 & \varepsilon \\
        \varepsilon & \varepsilon & \varepsilon & \varepsilon & \varepsilon & 1 & 4 & \varepsilon & \varepsilon &  1 \\
        \varepsilon & \varepsilon & \varepsilon & \varepsilon & \varepsilon & \varepsilon & \varepsilon & \varepsilon & 1 & \varepsilon 
    \end{pmatrix}.
\end{align*}
 \begin{figure}[tbp]
 \centering
 \label{fig:exm}\includegraphics[width=0.6\textwidth]{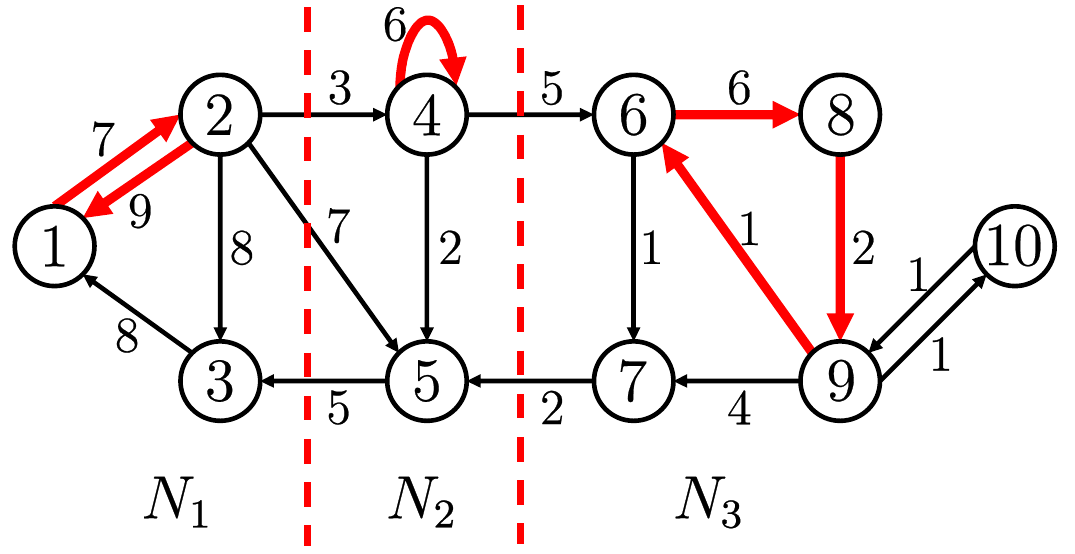}
 \caption{Graph $\mathcal{G}(A)$ in Section 4.4. Red-bold arrows represent quasi-critical circuits.}
 \end{figure}
    The graph $\mathcal{G}(A)$ is illustrated in Figure \ref{fig:exm}. The roots of the characteristic polynomial of $A$ and the corresponding MMCS are
    \begin{align*}
        \lambda_{1} = 8, \lambda_{2} = 7, \lambda_{3} = 6, \lambda_{4} = 4, \lambda_{5} = 3, \lambda_{6} = 0, \lambda_{7} = \varepsilon,
    \end{align*}
    and 
    \begin{align*}
        \mathfrak{C}_{0} &= \emptyset, \\
        \mathfrak{C}_{1} &= \{(1,2,1)\}, \\
        \mathfrak{C}_{2} &= \{(1,2,3,1)\}, \\
        \mathfrak{C}_{3} &= \{(1,2,3,1), (4,4)\}, \\
        \mathfrak{C}_{4} &= \{(1,2,5,3,1), (4,4)\}, \\
        \mathfrak{C}_{5} &= \{(1,2,5,3,1), (4,4), (6,8,9,6)\}, \\
        \mathfrak{C}_{6} &= \{(1,2,4,6,8,9,7,5,3,1)\}, 
    \end{align*}
    respectively. The partition of $N$ by Algorithm \ref{alg:partition} is 
    \begin{align*}
        N_{1} = \{1,2,3\},\ N_{2} = \{4,5\},\ N_{3} = \{6,7,8,9,10\}
    \end{align*}
    with $\lambda_{k(1)} = 8, \lambda_{k(2)} = 6$ and $\lambda_{k(3)} = 3$. The corresponding quasi-critical circuits are
    \begin{align*}
        \mathcal{C}_{1} = (1,2,1),\ \mathcal{C}_{2} = (4,4),\ \mathcal{C}_{3} = (6,8,9,6),
    \end{align*}
    respectively. Next, we apply Algorithm \ref{alg:visualize}. We start with $N_{3}$. After a trivial step $i =10$, we proceed as follows:
    \begin{align*}
        i &= 9: & b(N') &= (-\lambda_{k(3)}) \otimes A(N') = \begin{pmatrix}
        \varepsilon & -2 \\
        -2 & \varepsilon 
        \end{pmatrix}, \quad
        & w(N') &= \begin{pmatrix} 0 \\ -2 \end{pmatrix} \\ 
        i &= 8: & b(N') &= \begin{pmatrix}
        \varepsilon & -1 & \varepsilon \\
        \varepsilon & \varepsilon & -4 \\
        \varepsilon & 0 & \varepsilon 
        \end{pmatrix}, \quad
        & w(N') &= \begin{pmatrix} 0 \\ 0 \\ 0 \end{pmatrix} \\ 
        i &= 7: & b(N') &= \begin{pmatrix}
        \varepsilon & \varepsilon & \varepsilon & \varepsilon \\
        \varepsilon & \varepsilon & -1 & \varepsilon \\
        1 & \varepsilon & \varepsilon & -4 \\
        \varepsilon & \varepsilon & 0 &  \varepsilon 
        \end{pmatrix},
        & w(N') &= \begin{pmatrix} 0 \\ 0 \\ 1 \\ 1 \end{pmatrix}, \\
        i &= 6: & b(N') &= \begin{pmatrix}
        \varepsilon & -2 & 3 & \varepsilon & \varepsilon \\
        \varepsilon & \varepsilon & \varepsilon & \varepsilon & \varepsilon \\
        \varepsilon & \varepsilon & \varepsilon & 0 & \varepsilon \\
        -3 & 0 & \varepsilon & \varepsilon & -4 \\
        \varepsilon & \varepsilon & \varepsilon & 0 & \varepsilon 
        \end{pmatrix},
        & w(N') &= \begin{pmatrix} 0 \\ -3 \\ -3 \\ -3 \\ -3 \end{pmatrix}. 
    \end{align*}
    Here, $b(N')$ is a matrix consisting of $b((j_{1},j_{2}))$ for $j_{1},j_{2} \in N'$ just before the maximum weight path tree for $\mathcal{G}(N')$ is derived, and $w(N')$ is a vector consisting of $w_{j}$ for $j \in N'$ at the end of the iteration for $i$. Thus, we obtain    
    \begin{align*}
        A'_{3} = \begin{pmatrix}
        \varepsilon & -5 & 0 & \varepsilon & \varepsilon \\
        \varepsilon & \varepsilon & \varepsilon & \varepsilon & \varepsilon \\
        \varepsilon & \varepsilon & \varepsilon & 0 & \varepsilon \\
         0 & 0 & \varepsilon & \varepsilon & -4 \\
        \varepsilon & \varepsilon & \varepsilon & 0 & \varepsilon 
        \end{pmatrix}, \quad
         \bm{d}_{3} &= \begin{pmatrix} 0 \\ -3 \\ -3 \\ -2 \\ -4 \end{pmatrix}. 
    \end{align*}
    When moving to $N_{2}$, we set $b(N_{3}) = (\lambda_{k(3)} - \lambda_{k(2)}) \otimes A'_{3} = (-3) \otimes A'_{3}$.
    Continuing the computation, we obtain the visualizations
    \begin{align*}
        A'_{2} &= \begin{pmatrix}
        0 & -4 & -10 & \varepsilon & \varepsilon & \varepsilon & \varepsilon \\
        \varepsilon & \varepsilon & \varepsilon & \varepsilon & \varepsilon & \varepsilon & \varepsilon \\
        \varepsilon & \varepsilon & \varepsilon & 0 & -1 & \varepsilon & \varepsilon \\
        \varepsilon & 0 & \varepsilon & \varepsilon & \varepsilon & \varepsilon & \varepsilon \\
        \varepsilon & \varepsilon & \varepsilon & \varepsilon & \varepsilon & 0 & \varepsilon \\
        \varepsilon & \varepsilon &  -8 & 0 & \varepsilon & \varepsilon & -10 \\
        \varepsilon & \varepsilon & \varepsilon & \varepsilon & \varepsilon & 0 & \varepsilon 
        \end{pmatrix}, 
        & \bm{d}_{2} &= \begin{pmatrix} 0 \\ 0 \\ -9 \\ -4 \\ -10 \\ -6 \\ -11 \end{pmatrix}, \\        
        A'_{1} &= \begin{pmatrix}
        \varepsilon & 0 & \varepsilon & \varepsilon & \varepsilon & \varepsilon & \varepsilon & \varepsilon & \varepsilon & \varepsilon \\
        0 & \varepsilon & -1 & -15 & -5 & \varepsilon & \varepsilon & \varepsilon & \varepsilon & \varepsilon \\
        0 & \varepsilon & \varepsilon & \varepsilon & \varepsilon & \varepsilon & \varepsilon & \varepsilon & \varepsilon & \varepsilon \\
        \varepsilon &  \varepsilon & \varepsilon & -2 & 0 & -10 & \varepsilon & \varepsilon & \varepsilon & \varepsilon \\
        \varepsilon & \varepsilon & 0 & \varepsilon & \varepsilon & \varepsilon & \varepsilon & \varepsilon & \varepsilon & \varepsilon \\
        \varepsilon & \varepsilon & \varepsilon & \varepsilon & \varepsilon & \varepsilon & 0 & -5 & \varepsilon & \varepsilon \\
        \varepsilon & \varepsilon & \varepsilon & \varepsilon & 0 & \varepsilon & \varepsilon & \varepsilon & \varepsilon & \varepsilon \\
        \varepsilon & \varepsilon & \varepsilon & \varepsilon & \varepsilon & \varepsilon & \varepsilon & \varepsilon &  0 & \varepsilon \\
        \varepsilon & \varepsilon & \varepsilon & \varepsilon & \varepsilon &  -10 &  0 & \varepsilon & \varepsilon & -14 \\
        \varepsilon & \varepsilon & \varepsilon & \varepsilon & \varepsilon & \varepsilon & \varepsilon & \varepsilon & 0 & \varepsilon 
        \end{pmatrix}, 
        & \bm{d}_{3} &= \begin{pmatrix} 0 \\ 1 \\ 0 \\ -9 \\ -3 \\ -16 \\ -9 \\ -19 \\ -13 \\ -20 \end{pmatrix}.
    \end{align*}
    Subsequently, we compute matrices $C_{1},C_{2},C_{3}$ and $R_{1},R_{2},R_{3}$ for the CSR expansion. For the visualized matrix $A'_{3}$, we construct an extended graph $\hat{G}$ with $15$ nodes (shown in Figure~\ref{fig:exm2}). We find the maximum weight paths from $(6,0)$ to all nodes. The maximum weight to node $(i,k)$ is denoted by $d_{i,k}$. We obtain
    \begin{align*}
        \hat{R}_{3} = \begin{pmatrix} d_{6,0} & d_{7,0} & d_{8,0} & d_{9,0} & d_{10,0} \\ 
        d_{6,1} & d_{7,1} & d_{8,1} & d_{9,1} & d_{10,1} \\ 
        d_{6,2} & d_{7,2} & d_{8,2} & d_{9,2} & d_{10,2} \end{pmatrix}
        = \begin{pmatrix} 0 & 0 & -4 & -8 & -4 \\ 
        -8 & -5 & 0 & -4 & -12 \\ 
        -4 & -4 & -8 & 0 & -8 \end{pmatrix}
    \end{align*}
    and subsequently
    \begin{align*}
        R_{3} = \begin{pmatrix} \mathcal{E} & \hat{R}_{3} \otimes \mathrm{diag}(-\bm{d}_{3}) \end{pmatrix}= 
        \begin{pmatrix} \varepsilon & \varepsilon & \varepsilon & \varepsilon & \varepsilon & 0 & 3 &  -1 & -6 & 0 \\ 
        \varepsilon & \varepsilon & \varepsilon & \varepsilon & \varepsilon & -8 & -2 & 3 & -2 & -8 \\ 
        \varepsilon & \varepsilon & \varepsilon & \varepsilon & \varepsilon & -4 & -1 & -5 & 2 & -4 \end{pmatrix}.
    \end{align*}
    Similarly, we obtain
    \begin{align*}
        C_{3} &= \begin{pmatrix} \mathcal{E} & \mathrm{diag}(\bm{d}_{3}) \otimes \begin{pmatrix} 0 & \varepsilon & -8 & -4 & -8  \\ 
        -4 & \varepsilon & 0 & -8 & 0 \\ 
        -8 & \varepsilon & -4 & 0 & -4 \end{pmatrix} \end{pmatrix}^{\top} \\
        &= \begin{pmatrix} \varepsilon & \varepsilon & \varepsilon & \varepsilon & \varepsilon & 0 & \varepsilon & -11 & -6 & -12 \\ 
        \varepsilon & \varepsilon & \varepsilon & \varepsilon & \varepsilon & -4 & \varepsilon & -3 & -10 & -4 \\ 
        \varepsilon & \varepsilon & \varepsilon & \varepsilon & \varepsilon & -8 & \varepsilon & -7 & -2 & -8 \end{pmatrix}^{\top} \\
        R_{2} &= \begin{pmatrix} \mathcal{E} & \begin{pmatrix} 0 & -4 & -10 & -10 & -11 & -11 & -21 \end{pmatrix} \otimes \mathrm{diag}(-\bm{d}_{2}) \end{pmatrix}\\
        &= \begin{pmatrix} \varepsilon & \varepsilon & \varepsilon & 0 & -4 & -1 & -6 & -1 & -5 & -10 \end{pmatrix},\\
        C_{2} &= \begin{pmatrix} \mathcal{E} & \mathrm{diag}(\bm{d}_{2}) \otimes \begin{pmatrix} 0 & \varepsilon & \varepsilon & \varepsilon & \varepsilon & \varepsilon & \varepsilon \end{pmatrix} \end{pmatrix}^{\top} \\
        &= \begin{pmatrix} \varepsilon & \varepsilon & \varepsilon & 0 & \varepsilon & \varepsilon & \varepsilon & \varepsilon & \varepsilon & \varepsilon \end{pmatrix}^{\top},\\
        R_{1} &= \begin{pmatrix} 0 & -1 & -1 & -15 & -5 & -26 & -25 & -30 & -31 & -44 \\
        -1 & 0 & -1 & -16 & -6 & -25 & -26 & -31 & -30 & -45 \end{pmatrix} \otimes \mathrm{diag}(-\bm{d}_{1}) \\
        &= \begin{pmatrix} 0 & -2 & -1 & -6 & -2 & -10 & -16 & -11 & -18 & -24 \\
        -1 & -1 & -1 & -7 & -3 & -9 & -17 & -12 & -17 & -25 \end{pmatrix}, \\
        C_{1} &= \mathrm{diag}(\bm{d}_{1}) \otimes \begin{pmatrix} 0 & -1 & -1 & -1 & 0 & 0 & -1 & -1 & 0 & -1 \\
        -1 & -2 & 0 & 0 & -1 & -1 & 0 & 0 & -1 & 0 \end{pmatrix}^{\top} \\
        &= \begin{pmatrix} 0 & 0 & -1 & -10 & -3 & -16 & -10 & -20 & -14 & -21 \\
        -1 & -1 & 0 & -9 & -4 & -17 & -9 & -19 & -13 & -20 \end{pmatrix}^{\top}.
    \end{align*}
    Thus, the CSR expansion of $A$ is 
    \begin{align*}
        A^{\otimes t} &= 8^{\otimes t} \otimes C_{1} \otimes \begin{pmatrix} \varepsilon & 0 \\ 0 & \varepsilon \end{pmatrix}^{\otimes t} \otimes R_{1} \\
        &\qquad\qquad \oplus 6^{\otimes t} \otimes C_{2} \otimes R_{2} \oplus 3^{\otimes t} \otimes C_{3} \otimes \begin{pmatrix} \varepsilon & 0 & \varepsilon \\ \varepsilon & \varepsilon & 0 \\
        0 & \varepsilon & \varepsilon \end{pmatrix}^{\otimes t} \otimes R_{3}.
    \end{align*}

\section{Concluding remarks}

In this study, we proposed a new scheme for the CSR expansion that induces an $O(n(m+n\log n))$ algorithm for an $n \times n$ max-plus matrix, where $m$ is the number of finite entries. We showed that the growth rate of each term in the expansion is a root of the characteristic polynomial of the matrix. Thus, it suggests a new application of the characteristic polynomial to max-plus matrix theory. As described in~\cite{Nishida2020}, we can define algebraic eigenvectors corresponding to the roots of the characteristic polynomial. We expect that these vectors generalizes the CSR expansion by providing some information on the dynamics before becoming periodic, which remains as a future work.

\section*{Acknowledgments}
This work was funded by JSPS KAKENHI No.22K13964.

\bibliographystyle{abbrv}
\bibliography{CSRexpand_arXiv}
%\begin{thebibliography}{99}
%\bibitem{b1} Bib1
%\end{thebibliography}

\end{document}